\documentclass[12pt]{article}
\usepackage{amsthm}
\usepackage{amssymb}
\usepackage{graphicx}
\usepackage{times}
\usepackage[active]{srcltx}
\newcommand{\indep}{\;\, \rule[0em]{.03em}{.67em} \hspace{-.25em}
\rule[0em]{.65em}{.03em} \hspace{-.25em}
\rule[0em]{.03em}{.67em}\;\,}

\newcommand{\heading}[1]{\bigskip{\noindent\large\bf #1}}

\newcommand{\tr}{{\rm tr}}

\newcommand{\E}{{\rm E}}

\newcommand{\T}{^{\rm T}}

\newcommand{\eye}{{\mathcal{I}}}

\newcommand{\real}{\mathbb{R}}

\newcommand{\f}{{l}}
\renewcommand{\k}{{f}}

\newtheorem{theorem}{\sc Theorem}
\newtheorem{lemma}{\sc Lemma}

\title{The Lineartiy Condition and Adaptive Estimation in Single-index Regressions}
\author{Yongwu Shao\footnote{Corresponding author.  E-mail: ywshao@stat.umn.edu} \footnote{R. D. Cook and Y. Shao were suppoted by Grant DMS-0405360 from the National Science Foundation.}\and Dennis Cook\footnotemark[\value{footnote}]\and Sanford Weisberg\\
{\small\it School of Statistics, University of Minnesota, Minneapolis, MN 55455, USA.}}
\date{Sep 15, 2006}

\begin{document}

\maketitle

\begin{abstract}
We show that under a linearity condition on the distribution of the predictors, the coefficient vector in a single-index regression can be estimated with the same efficiency as in the case when the link function is known.  Thus, the linearity condition seems to substitute for knowing the exact conditional distribution of the response given the linear combinations of the predictors.
\end{abstract}

\section{Introduction}
\subsection{Single-index Regressions}
Consider a continuous univariate response $Y$ and a vector of continuous predictors $X\in\real^p$.  The most general goal of a regression is to infer about the conditional distribution of $Y|X$.  In this paper we consider single-index regressions, in which $Y|X$ depends on $X$ through at most one linear combination $\beta_0\T X$ of the predictors.

Focusing on the mean function $\E(Y|X)$, H\"ardle and Stoker (1989) developed a nonparametric method called average derivative estimation for estimating $\beta_0$ in the single-index conditional mean $\E(Y|X)=g(\beta_0\T X)$, where the mean function $g$ is unknown.  Weisberg and Welsh (1994) considered the case in which $Y|X$ follows a generalized linear model, where the linear coefficient $\beta_0$ and the link function are unknown.  Both pairs of authors gave estimates for $\beta_0$ that are $\sqrt{n}$-consistent.

Yin and Cook (2005) proposed the problem of single-index regressions, in which the conditional distribution of $Y|X$ is completely characterized by a linear combination $\beta_0\T X$, so there is no loss of information about $Y$ if we replace $X$ with $\beta_0\T X$.  More specifically, we assume that
\begin{equation}\label{def:cs}
Y\indep X|\beta_0\T X
\end{equation}
where for identifiability purposes, we require that $||\beta_0||=1$.  (\ref{def:cs}) is equivalent to the statement that $Y|X$ has a conditional density $\eta_0(y|\beta_0\T x)$, where $\eta_0$ is unknown.  Single-index regression is a special case of sufficient dimension reduction when the dimension of the central subspace (Cook, 1996) is one.  It does not require a pre-specified single-index model.

We show that under the linearity condition (Li and Duan, 1989), for single-index regressions there exists an adpative estimate for $\beta_0$ that can be estimated with the same efficiency as the maximum likelihood estimate when the conditional density $\eta_0$ is completely specified.  For example, if the true model is $Y=g(\beta_0\T X)+\epsilon$, where the link function $g$ and the density of the error $\epsilon$ are unknown, then $\beta_0$ can be estimated with the same efficiency as in the case when $g$ and the error distribution are known.

\subsection{Linearity Condition}
Many sufficient dimension reduction methods require the linearity condition: $E(X|\beta_0\T X)$ is a linear function of $\beta_0\T X$ (Li and Duan, 1989).  It is used in popular methods like sliced inverse regression (Li, 1991), sliced average variance estimation (Cook and Weisberg, 1991) and principal Hessian directions (Li, 1992, Cook, 1998).  

The linearity condition holds if the predictor has an elliptically distribution (Eaton, 1986), so it holds when $X$ has a multivariate normal distribution.  Moreover, Diaconis and Freedman (1984) showed that most low-dimension projections of a high-dimension data cloud are close to being normal.  Hall and Li (1993) argued that the linearity condition holds approximately when $p$ is large.  The linearity condition applies only to the marginal distribution of the predictors and not to the conditional distribution of $Y|X$ as is common in regression modeling.  Consequently at the stage of data collection, we might design the experiment so that the distribution of $X$ will not blatantly violate elliptic symmetry.  We can also transform the predictors to normality, or we can re-weight the data (Cook and Nachtscheim, 1994) to approximate an ellliptical distribution.

\subsection{Adaptive Estimation}

The problem of adaptive estimation was introduced by Stein (1956).  One wishes to estimate a Euclidean parameter $\theta$ in the presence of an infinite-dimensional shape parameter $G$, usually the density.  An adaptive estimate performs asymptotically as well with $G$ unknown as the maximum likelihood estimate does when $G$ is known (Bickel, 1982).  A general method of constructing adaptive estimates was constructed by Bickel (1982).  Schick (1986, 1993) generalized and improved Bickel's method.

It has been shown that adaptive estimation is possible in the symmetric location problem, in which we need to estimate the center of symmetry of an unknown distribution (Stone, 1975).  It is also possible in linear regressions where the error density is symmetric and unknown and we need to estimate the linear coefficient (Bickel, 1982).  When the observations are not independent, Koul and Pflug (1990), Schick (1993), Koul and Schick (1996) showed adaptive estimation is possible in certain autoregressive models.  Early literature on adaptive estimation generally focused on these models and their generalizations.  In this paper we show that under the linearity condition, adaptive estimation is also possible for single-index regressions.

\section{Main Results}
Without loss of generality we assume that $X$ has mean zero and covariance $I_p$.  We also assume that $\beta_0\in\Theta$, where $$\Theta=\{\beta\in\real^p:||\beta||=1\}.$$
Let $\f(t,y)=(\partial/\partial t)\eta_0(y|t)/\eta_0(y|t)$ be the derivative of the log density or equivalently the log likelihood.  By using a Lagrange multiplier, the score equation for $\beta_0$ is 
\begin{equation}\label{eq:TrueScore}
	Q_{\beta_0}\E[X\f(\beta_0\T X,Y)]=0,
\end{equation}
where $Q_{\zeta}=I_p-P_{\zeta}$ and $P_\zeta$ is the orthogonal projection onto the subspace spanned by the columns of the matrix $\zeta$.

It can be shown that (\ref{eq:TrueScore}) holds not only for $\f$, but for any $\k(\cdot,\cdot)\in\real$.

\begin{lemma}\label{lemma:1}
Assume that the linearity condition holds.  Assume $\k(\cdot,\cdot)\in\real$.  Then $\beta_0$ is a solution of the equation \begin{equation}\label{eq:unique} Q_\beta\E[X\k(\beta\T X,Y)]=0. \end{equation}
\end{lemma}

\begin{proof}
Since $X$ has covariance matrix $I_p$, according to Cook (1998, pp. 57), we have $\E[Q_{\beta_0}X|\beta_0\T X]=0$.  Therefore
\begin{eqnarray*}
Q_{\beta_0}\E[X\k(\beta_0\T X,Y)] &=& \E[Q_{\beta_0}X\k(\beta_0\T X,Y)] \\
&=& \E\{\E[Q_{\beta_0}X\k(\beta_0\T X,Y)]|\beta_0\T X \}\\
&=& \E\{\E[Q_{\beta_0}X|\beta_0\T X]\E[\k(\beta_0\T X,Y)|\beta_0\T X]\} \\
&=& 0
\end{eqnarray*}
\end{proof}

The above lemma shows that a misspecified $\f$ still produces a Fisher consistent estimate of $\beta_0$.  According to van der Vaart (1998, Theorem 25.27), Lemma~\ref{lemma:1} together with some regularity conditions would enable us to construct an adaptive estimate for $\beta_0$.  The regularity conditions are typically satisfied in practice.  A proof of the following theorem is given in the appendix.

\begin{theorem}\label{thm:1}
Assume that the Fisher information $\eye(\beta)=\E_\beta[XX\T \f^2(\beta\T X,Y)]$ is finite, nonsingular and differentiable with respect to $\beta$ in a neighborhood of $\beta_0$.  Let $\hat\f_n(t,y)$ be an estimate of $\f(t,y)$ that satisfies 
\begin{equation}\label{eq:cond} \E_{\beta_0}[||X||^2(\hat\f(\beta_0\T X,Y)-\f(\beta_0\T X,Y))^2]=o_p(1). \end{equation}
Then under the linearity condition we can construct an adaptive estimate of $\beta_0$ in (\ref{def:cs}) based on $\hat\f_n(t,y)$.
\end{theorem}

Following van der Vaart (1998, pp. 393), an adaptive estimate can be constructed in the following way.  Suppose $\beta_n$ is a $\sqrt{n}$-consistent estimate of $\beta_0$.  For instance, under the linearity condition $\beta_n$ can be chosen as the ordinary least squares estimator (Li and Duan, 1989).  Let $\Gamma_n$ be a $p\times(p-1)$ matrix such that $(\Gamma_n,\beta_n)$ is an orthogonal matrix. 
Let \[ \tilde \eye_n=\sum_{i=1}^n[X_iX_i\T\hat\f_n^2(\beta_n\T X_i,Y_i)] \] be an estimator of the information matrix for $\beta$.  Let $\hat\beta_n$ be a one-step iteration of the Newton-Raphson algorithm for solving the equation \[ Q_\beta \sum_{i=1}^n [X_i\hat\f_n(\beta\T X_i,Y_i)]=0 \] with respect to $\beta$ on the manifold $\Theta$, starting at the initial guess $\beta_n$.  We can write $\hat\beta_n$ as
\begin{equation}\label{eq:hatBeta}
	\hat\beta_n=\beta_n+\frac{1}{n}\Gamma_n[\Gamma_n\T\tilde \eye_n\Gamma_n]^{-1}\Gamma_n\T\sum_{i=1}^n [X_i\hat\f_n(\beta_n\T X_i,Y_i)]
\end{equation}
Van der Vaart (1998, Theorem 25.27) showed that, by using discretization and sample-splitting devices, $\hat\beta_n$ is an adaptive estimate of $\beta_0$ if $\hat\f$ satisifies (\ref{eq:cond}).  One such $\hat\f$ based on the kernel density estimation in H\"ardle and Stoker (1989) is constructed in the Appendix.

Since $\hat\beta_n$ is an adaptive estimator, it has the same asymptotic distribution as the maximum likelihood estimator.  Next we will derive the asymptotic distribution of the maximum likelihood estimator.  Let $\hat\beta_{\rm mle}$ be the maximum likelihood estimator of $\beta_0$.  It is shown in the appendix that under mild regularity conditions $\hat\beta_{\rm mle}$ has the following asymptotic distribution.
\begin{theorem}\label{thm:2}
Assume that the regularity conditions for the asymptotic normality of the maximum likelihood estimate hold.  Then
\begin{equation}\label{eq:asym1}
\hat\beta_{\rm mle}=\beta_0+\frac{1}{n}\Gamma_0[\Gamma_0\T \eye(\beta_0)\Gamma_0]^{-1}\Gamma_0\T \sum_{i=1}^n[X_i\f(\beta_0\T X_i,Y_i)]+o_p(n^{-1/2})
\end{equation}
where $\Gamma_0$ is a $p\times(p-1)$ matrix such that $(\Gamma_0,\beta_0)$ is an orthogonal matrix.
\end{theorem}
Since $\hat\beta_n$ is an adaptive estimator, it has the same asymptotic distribution as $\hat\beta_{\rm mle}$, we conclude that $\sqrt{n}(\hat\beta_n-\beta_0)$ converges to a normal distribution with zero mean and covariance matrix equal to the covariance matrix of $\Gamma_0[\Gamma_0\T \eye(\beta_0)\Gamma_0]^{-1}\Gamma_0\T X\f(\beta_0\T X,Y)]$.

\section{Discussion}

In this article we showed that under the linearity condition, there exists an adaptive estimate of the coefficient vector in a single-index regression.  From this result we can see the important role of the lineartiy condition in single-index regression, and more generally, in sufficient dimension reduction.  The linearity condition is unusual, as it does not occur commonly outside of sufficient dimension reduction.  We have shown that the linearity condition asymptotically takes the place of a known density.  We conjecture that if the linearity condition fails, then an adaptive estimate does not exist.  As a consequence, the coefficient vector cannot be estimated as well as it can be with the maximum likelihood estimator.


\newpage
\heading{Appendix}

\begin{proof}[Proof of Theorem~\ref{thm:1}]
Since Lemma~\ref{lemma:1} holds, according to van der Vaart (1998, Theorem 25.27), we only need to prove the following two statements.
\begin{enumerate}
\item The conditional density $\eta_0(y|\beta_0\T x)$ is differentiable in quadratic mean with respect to $\beta_0$.
\item Let $h(\beta\T x,y)$ be the joint density of $\beta\T X$ and $Y$, then \[ \int ||x||^2\left[\f(\beta_n\T x,y)\sqrt{h(\beta_n\T x,y)}-\f(\beta_0\T x,y)\sqrt{h(\beta_0\T x,y)}\right]^2dxdy\rightarrow 0. \]
\end{enumerate}

The first statement is true by van der Vaart (1998, Theorem~7.2).  So we only need to prove the second statement.

Since
\[ \eye(\beta_0)=\int xx\T\left[\f(\beta_0\T x,y)\sqrt{h(\beta_0\T x,y)}\right]^2dxdy \]
and
\[ \eye(\beta_n)=\int xx\T\left[\f(\beta_n\T x,y)\sqrt{h(\beta_0\T x,y)}\right]^2dxdy \]
By the assumptions, $\eye(\beta)$ is continuous on a neighborhood of $\beta_0$ and $\eye(\beta_0)$ is finite, we conclude that $\eye(\beta_n)$ is also finite, hence $\tr[\eye(\beta_0)+\eye(\beta_n)]<\infty$.  By the triangular inequality,
\begin{eqnarray*}
&& \int ||x||^2\left[|\f(\beta_n\T x,y)|\sqrt{h(\beta_n\T x,y)}+|\f(\beta_0\T x,y)|\sqrt{h(\beta_0\T x,y)}\right]^2dxdy \\
&\leq&\tr[\eye(\beta_0)+\eye(\beta_n)]<\infty 
\end{eqnarray*}
Then by the dominate convergence theorem, 
\[ \int ||x||^2\left[\f(\beta_n\T x,y)\sqrt{h(\beta_n\T x,y)}-\f(\beta_0\T x,y)\sqrt{h(\beta_0\T 
x,y)}\right]^2dxdy\rightarrow 0 \]
Therefore the second statement is also true.

\end{proof}

\begin{proof}[Proof. of Theorem~\ref{thm:2}]
We first transform  the manifold $\Theta$ to $\real^{p-1}$ by using the following linear transformation.  For any $\beta\in\Theta$, let $\alpha=\varphi(\beta)=\Gamma_0\beta$.  Then $\beta=\varphi^{-1}(\alpha)=\Gamma_0\alpha+(1-||\alpha||^2)\beta_0$, and $\eta_0(y|\beta\T x)=\eta_0(y|\varphi^{-1}(\alpha)\T x)$.  By taking the derivative of $\eta_0(y|\varphi^{-1}(\alpha)\T x)$ with respect to $\alpha$, we can derive the asymptotic distribution of the maximum likelihood estimate for $\alpha$ as following,
\[ \hat\alpha_{\rm mle}=\frac{1}{n}[\Gamma_0\T \eye(\beta_0)]^{-1}\Gamma_0\T \sum_{i=1}^n[X_i\f(\beta_0\T X_i,Y_i)]+o_p(n^{-1/2}) \]
By the delta method, we have
\[ \hat\beta_{\rm mle}=\beta_0+\frac{1}{n}\Gamma_0[\Gamma_0\T \eye(\beta_0)\Gamma_0]^{-1}\Gamma_0\T \sum_{i=1}^n[X_i\f(\beta_0\T X_i,Y_i)]+o_p(n^{-1/2}) \]

\end{proof}

\begin{proof}[Construction of $\hat\f$ that satisfies (\ref{eq:cond}).]
Let $h(\beta_0\T x,y)$ be the joint density of $(\beta_0\T X,Y)$, and $g(\beta_0\T x)$ be the density of $\beta_0\T X$, then $$\eta_0(y|\beta_0\T x)=h(\beta_0\T x,y)/g(\beta_0\T x)$$ and $$\f=h'/h-g'/g,$$  where $h',g'$ are the derivative of $h,g$ w.r.t. the first argument.  To estimate $\f$, we only need to estimate $h'/h$ and $g'/g$.  

We only consider the estimation of $g'/g$ in detail here, because $h'/h$ can be estimated in the same way, except that the dimension of the density estimation is different.  Let $d$ be the dimension of the density estimation, $d=1$ for $g$ and $d=2$ for $h$.

Let $T_i=\beta_0\T X_i$.  For a fixed twice continuously differentiable probability density $w$ with compact support, a bandwidth parameter $\sigma$, and a cut-off tuning parameter $\delta$, set
\begin{eqnarray}
\hat g_n(s) &=& \sigma_n^{-d}\sum_{i=1}^n w(\frac{s-T_i}{\sigma_n}) \nonumber\\ 
\hat \xi_n(s) &=& \frac{\hat g_n'}{\hat g_n}(s)1_{\hat g_n(s)>\delta} \label{eq:kn} 
\end{eqnarray}
where $\hat \xi_n(s)$ is our estimator of $g'(s)/g(s)$.  Then $\E[(\hat\xi_n(X)-g'(X)/g(X))^2||X||^2]$ converges to zero in probability provided $\delta\uparrow\infty$ and $\sigma\downarrow 0$ at appropriate speeds.  

Hardle and Stoker, 1991, page 992) showed that under some regularity conditions we have for any $\epsilon>0$,
\[ \sup[|\hat g(s)-g(s)|1_{g(s)>(\delta/2)}]=O_p[(n^{1-(\epsilon/2)}\sigma^d)^{-1/2}] \]
and
\[ \sup[|\hat g'(s)-g'(s)|1_{g(s)>(\delta/2)}]=O_p[(n^{1-(\epsilon/2)}\sigma^{d+2})^{-1/2}] \]
Therefore
\[ \sup[|(\hat g'/\hat g)-(g/g)|1_{g>(\delta/2)}]=O_p[\delta^{-2}(n^{1-(\epsilon/2)}\sigma^{d+2})^{-1/2}] \]
Hence for large $n$ we have
\begin{eqnarray*}
&&\E[(\hat\xi_n-g')^2||X||^2] \\
&=& \E[(g'/g)^2||X||^2 1_{\hat g<\delta}] +\E[((\hat g'/\hat g)-(g'/g))^2||X||^2 1_{\hat g>\delta}] \\
&\leq& \E[(g'/g)^2||X||^2 1_{g<2\delta}] +\E[((\hat g'/\hat g)-(g'/g))^2||X||^2 1_{g>(\delta)/2}] \\
&\leq& \E[(g'/g)^2||X||^2 1_{g<2\delta}] + O_p[\delta^{-2}(n^{1-(\epsilon/2)}\sigma^{d+2})^{-1/2}]\cdot\E[||X||^2]
\end{eqnarray*}
Assume that $\E[||X||^2]<\infty$.  Since $(g'/g)^2||X||^2 1_{g<2\delta}$ is dominated by $(g'/g)^2||X||^2$, and $\E[(g'/g)^2||X||^2 1_{g<2\delta}]$ is finite by assumptions, therefore $\E[(g'/g)^2||X||^2 1_{g<2\delta}]$ converges to zero when $\delta$ goes to zero.  By the assumptions, $\delta^{-2}(n^{1-(\epsilon/2)}\sigma^{d+2})^{-1/2}$ also converges to zero, therefore $\E[(\hat\xi_n-(g'/g))^2||X||^2]=o_p(1)$.

In the same fashion we can construct $\hat\zeta_n(s,y)$ to estimate $h'/h$, except that we use $(T_i,Y_i)$ as observations.  Then an estimator for $\f$ can be defined as \begin{equation}\label{eq:kn1} \hat\f_n=\hat \zeta_n-\hat\xi_n. \end{equation}  Since $\E[(\hat\xi_n-g'/g)^2||X||^2]$ and $\E[(\hat\zeta_n-h'/h)^2||X||^2]$ converges to zero in probability, we have $\E[(\hat\f_n(\beta_0\T X,Y)-f(\beta_0\T X,Y))^2||X||^2]$ converges to zero in probability, and (\ref{eq:cond}) is satisfied.
\end{proof}


\newcommand{\bibit}[5]
{\item[] {#1} #2. {#3}, {#4}, #5.}

\newcommand{\bibitbook}[3]
{\item[] {#1} {#2}. #3.}

\newcommand{\ann}{Ann. of Statist.}
\newcommand{\jasa}{J. Amer. Statist. Assoc.}
\newcommand{\jrssb}{J. Roy. Statist. Soc. Ser. B}

\heading{References}
\small
\begin{description}

\renewcommand{\&}{{\rm and}}

\bibit
{Bickel, P. J. (1982)} {On adaptive estimation} {\ann} {10} {647-671}

\bibit
{Cook, R. D. (1996)} {Graphics for regressions with a binary response} {\jasa} {91} {983-992}

\bibit
{Cook, R. D. (1998)} {Principle Hessian directions revisited (with discussion)} {\jasa} {93} {84-100}

\bibit
{Cook, R. D. \& Nachtsheim, C. J. (1994)} {Re-weighting to achieve elliptically contoured covariates in regression} {\jasa} {89} {592-600}

\bibit
{Cook, R. D. \& Weisberg, S. (1991)}  {Discussion of "sliced
inverse regression for dimension reduction"} {\jasa} {86} {328-332}

\bibit
{Diaconis, P. \& Freedman, D. (1984)} {Asymptotics of graphical projection pursuit} {\ann} {12} {793-815}

\bibit
{Eaton, M. L. (1986)} {A characterization of spherical distributions} {J. Mult. Anal.} {20} {272-6}

\bibit
{Hall, P. \& Li, K. C. (1993)} {On almost linearity of low dimensional projections from high dimensional data} {\ann} {21} {867-889}

\bibit
{H\"ardle, W. \& Stoker, T. M. (1989)} {Investigating smooth multiple regression by the method of average derivatives} {\jasa} {84} {986-995}

\bibit
{Koul, H. L. \& Pflug, G. (1990)} {Weakly adatptive estimators in explosive regression} {\ann} {18} {939-960}

\bibit
{Koul, H. L. \& Schick, A. (1996)} {Adaptive estimation in a random coefficient autoregressive model} {\ann} {24} {1025-1052}

\bibit
{Li, K. C. (1991)} {Sliced inverse regression for dimension
reduction (with discussion)} {\jasa} {86} {316-342}

\bibit
{Li, K. C. (1992)} {On Principle Hessian directions for data visualization and dimension reduction: Another application of Stein's lemma} {\jasa} {87} {1025-1039}

\bibit
{Li, K. C. \& Duan, N. (1989)} {Regression analysis under link violation} {\ann} {17} {1009-1052}

\bibit
{Schick, A. (1986)} {On asymptotically efficient estimation in semi-parametric models} {\ann} {14} {1139-1151}

\bibit
{Schick, A. (1993)} {On efficient estimation in regression models} {\ann} {21} {1481-1521}

\bibit
{Stein, C. (1956)} {Efficient nonparametric testing and estimation} {Proc. Third Berkeley Symp. Math. Statist. Prob.} {1} {187-196} Unversity of California Press.

\bibit
{Stone, C. J. (1975)} {Adaptive maximum likelihood estimators of a location parameter} {\ann} {3} {276-284}
\bibitbook
{van der Vaart, A. W. (1998)}{Asmyptotic Statistics}{Cambridge}

\bibit
{Weisberg, S. \& Welsh, A. H. (1994)} {Adapting for the missing link} {\ann} {22} {1674-1700}

\bibit
{Yin, X. \& Cook, R. D. (2005)} {Direction estimation in single-index regressions} {Biometrika} {92} {371-384}

\end{description}

\end{document}